# Monotonicity Property for a Class of Semilinear Partial Differential Equations.

Siva Athreya

March 5, 1999

**Abstract**

We establish a monotonicity property in the space variable for the solutions of an initial boundary value problem concerned with the parabolic partial differential equation connected with super-Brownian motion.

## 1  Introduction

The "hot spots" conjecture of J. Rauch has been analyzed in certain planar domains $D$ by R. Bañuelos and K. Burdzy [BB97] using probabilistic methods. In that paper, they synchronously couple two reflected Brownian motions and establish some monotonicity properties for solutions of the heat equation in $D$.

In this note we show that by applying similar coupling techniques to super-reflected Brownian motion one can prove a monotonicity property for solutions of a class of semilinear elliptic partial differential equations connected with super-Brownian motion. The result follows easily via an application of the existing machinery developed in the field of super-processes.

The purpose of this note is to enunciate the ease with which the probabilistic argument shown in [BB97] can be extended to provide a non-trivial result for solutions of certain semilinear partial differential equations. To the best of our knowledge the result presented in this note is new in the field of semilinear partial differential equations.

We consider solutions $u : \mathbb{R}_+ \times D \to \mathbb{R}_+$ of the following initial boundary value problem:

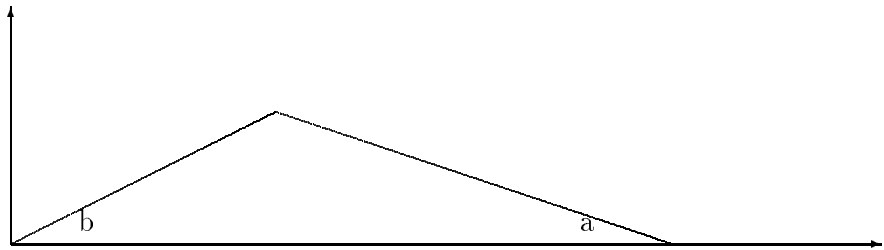

Figure 1: Obtuse triangle

$$\frac{\partial u}{\partial t}(t,x) = \frac{1}{2}\triangle u(t,x) + \Phi(u(t,x)), \quad x \in D, \quad t > 0, \tag{1}$$

$$u(0,x) = \phi(x), \quad x \in D, \tag{2}$$

$$\frac{\partial u}{\partial n}(t,x) = 0, \quad x \in \partial D, \quad t > 0. \tag{3}$$

Here $D$ is a bounded connected subset of $\mathbb{R}^2$, $\phi \in C^1$, and $\Phi : \mathbb{R}_+ \to \mathbb{R}_+$ is a continuously differentiable function of the form

$$\Phi(\lambda) = a_1 \lambda - b_1 \lambda^2 + \int_0^\infty (1 - \exp(-\lambda u) - \lambda u) \nu(du), \tag{4}$$

where $a_1 \in \mathbb{R}, b_1 \geq 0$ and $\nu(du)$ is a regular Borel measure in $\mathbb{R}_+$, such that $\int_0^\infty u \wedge u^2 \nu(du) < \infty$. It is well known that solutions to this initial-boundary problem exist and are unique.

The partial differential equations that arise when $\Phi(\lambda) = -\lambda^2$ (choose $b_1 = 1$ and $\nu(du) = 0$) and $\Phi(\lambda) = -\lambda^{1+\beta}, 0 < \beta < 1$ (choose $b_1 = 0$ and $\nu(du) = c_1 du/u^{2+\beta}$, $c_1 > 0$) are connected with binary and $\beta$-branching super-Brownian motion. The path properties of the process and various analytical properties of the partial differential equations have been extensively studied ([DIP89], [LG95] [Dyn91]).

Our main result is concerned with the direction of the gradient of $u(t, x)$ in obtuse triangles. We consider an obtuse triangle $D$, with the longest side of the triangle lying on the horizontal axis. The triangle lies in the first quadrant and one of its vertices is at the origin. The smaller sides of the triangle form angles $a$ and $b$ with the horizontal axis, with $a \in (-\frac{\pi}{2}, 0)$ and $b \in (0, \frac{\pi}{2})$ (See Figure 1). Let $\angle \nabla_x u(t, x)$ be the angle formed by the gradient $\nabla_x u(t, x)$ with the horizontal axis.

**Theorem 1.1.** *Suppose that $u(0, x)$ is $C^1$ and $c < \angle \nabla_x u(0, x) < d$ for all $x \in D$, where $c > b - \frac{\pi}{2}$ and $d < \frac{\pi}{2} + a$. Then for every $t$ and $x$ we have*

$$min(a, c) \leq \angle \nabla_x u(t, x) \leq max(b, d).$$

The main idea of the proof is to construct a synchronous coupling of historical reflected Brownian motions in $D$. We use the same method as in [Kle89]. The final step uses the log-Laplace functional of super-reflected Brownian motion to obtain the monotonicity property for the solutions to (1).

**Notation:** We shall denote $x \in \mathbb{R}^2$ as $x = (x^1, x^2)$, where each $x^i \in \mathbb{R}$ (real numbers). For any Polish space $G$, $y \in G$, measurable function $f : G \to R$ and a measure $m$ on $G$, we define $\langle f, m \rangle = \int_G f(y) dm(y)$ and $\delta_y$ as the dirac measure at the point $y$. $\mathcal{B}_d$ will denote the Borel $\sigma$-field on $\mathbb{R}^d$, $C_d = C([0,\infty), \mathbb{R}^d)$, $C_d^t = \{y \in C_d : y = y(\cdot \wedge t)\}$, $\mathcal{C}_d$ will denote the Borel $\sigma$-field of $C_d$, $M_F(C_d)$ the set of all finite measures on $C_d$, and $M_F(C_d^t) = \{m \in M_F(C_d) : y(\cdot \wedge t) = y \ m \ \text{a.e.} \ y\}$. For $z, w \in C_d$, we define

$$(z/s/w)(u) = \begin{cases} z(u) & \text{if } u < s, \\ w(u-s) & \text{if } u \geq s. \end{cases}$$

# 2 Synchronous coupling of historical reflected Brownian motions

First we provide a brief construction of "synchronous coupling" of reflected Brownian motions. We refer the reader to [BB97] for further details.

Let $B_t = (B_t^1, B_t^2)$ be a 2 dimensional Brownian motion starting at $x = (x^1, x^2)$, where $x^2 > 0$ and $C_t = (C_t^1, C_t^2) = (B_t^1 + (y^1 - x^1), B_t^2 + (y^2 - x^2))$, where $y = (y^1, y^2)$ with $y^2 > 0$. Define $\xi_t^x = (B_t^1, B_t^2 - 0 \wedge \min_{s \leq t} B_s^2)$ and $\xi_t^y = (C_t^1, C_t^2 - 0 \wedge \min_{s \leq t} C_s^2)$.

We shall call the pair $(\xi_t^x, \xi_t^y)$ a synchronous coupling of reflected Brownian motions in the upper half-plane. The above construction can be generalized to any polygonal domain $D \subset \mathbb{R}^2$. The construction gives us for every pair of $x, y \in D$, a pair of reflected Brownian motions $(\xi_t^x, \xi_t^y)$ starting at $(x, y)$, such that $\xi_t^x - \xi_t^y$ remains constant in the time periods when both processes are in the interior of the domain.

In this section, we shall use the synchronous coupling of two reflected Brownian motions in $D$ mentioned above as the path process to construct a coupling of two historical reflected Brownian motions.

We proceed to define a historical process $H_t$ on $D \times D$. Let $\xi_t^x, \xi_t^y$ be spatially coupled reflected Brownian motions as in [BB97], starting at $x$ and $y$ in $D$. Then $\xi_t = (\xi_t^x, \xi_t^y)$ is a continuous Markov process taking values in $D \times D$. The path valued process $\bar{\xi}_t = (\xi^x(\cdot \wedge t), \xi^y(\cdot \wedge t))$ will be the **motion process** for $H_t$. The **branching mechanism** is given by $\Phi$ described earlier in (4).

Applying Theorem 2.2.3. in [DP91] to the process $(\bar{\xi}, \Phi)$, we see that the $(\xi, \Phi)$-historical process $H_t$ exists. The semi-group $Q_{r,t}$ of $H$ is determined by

$$Q_{r,t}(\exp(-\langle m, \psi\rangle)) = \exp(-\langle m, V_{r,t}(\psi)\rangle),$$

where $V_{r,t}(\psi(z)) = \langle \delta_z, V_{r,t}(\psi)\rangle$ is the unique solution of

$$V_{r,t}(\psi(z)) = P_{r,y(r)}(\psi(z/r/\bar{\xi}_{t-r})) + \int_0^{t-r} P_{r,y(r)}(\Phi(V_{s+r,t}(z/r/\bar{\xi}_s)))ds, \tag{5}$$

where $\psi$ is a bounded positive measurable function, $m \in M_F(C_4)^r$ and $z \in C_4^r$, $r \leq t$.

Since the branching mechanism is spatially homogeneous, the measures $H_t^1(dx) = H_t(dx \times C_2)$ and $H_t^2(dy) = H_t(C_2 \times dy)$ are the historical processes associated with $\xi_t^x$ and $\xi_t^y$. We call the pair $\{H_t^1, H_t^2\}$ as synchronously coupled historical Brownian motion.

Let $\Pi_t(y) = y(t)$ be the coordinate map in $C_2$. Keeping the particle picture in mind, $X_t = H_t^1 \Pi_t^{-1}$ and $Y_t = H_t^2 \Pi_t^{-1}$ are 2-dimensional super-reflected Brownian motions in $\mathbb{R}^2$ with the branching mechanism $\Phi$ and the motion process $\xi^x$ and $\xi^y$ respectively. This is verified in Theorem 2.2.4 [DP91].

Let $\phi$ be a positive bounded measurable function in $D$, log-Laplace functional of $X_t$ is given by

$$\langle m, V_t\phi\rangle = -\log Q_m(\exp(-\langle X_t, \phi\rangle)), \quad m \in M_F(D), \tag{6}$$

where $v(t,x) = \langle \delta_x, V_t\phi\rangle = V_t\phi(x)$ is the the unique solution of

$$v(t,x) = S_t(\phi(x)) + \int_0^t S_t(\Phi(v(t-u,\cdot))(x)du; \tag{7}$$

Details can be found in [Fit88]. As $\phi$ and $\Phi$ are continuously differentiable, by Theorem 1.5 (page 187) in [Paz83], we conclude that $v(t,x)$ is a solution to (1), with initial conditions (2). As $v$ solves (7) and the fact that $S_t$ is the semi-group of reflected Brownian motion imply that $v$ satisfies the Neumann boundary conditions (3).

**Proof of Theorem 1.1:** Let $H_t^1$ and $H_t^2$ be the coupled historical reflected Brownian motion described above. Take any $(w,z)$ in the support of the historical measure $H_t$ and let $K(s)$ be the line joining the points $w(s)$ and $z(s)$.

Let $x, y$ in $D$ be chosen so that they satisfy

1. $x^1 < y^1$,

2. they are on a line $K$ inside the domain $D$, such that

$$\max(b,d) - \frac{\pi}{2} \leq \angle K \leq \min(a,c) + \frac{\pi}{2}. \tag{8}$$

Consider $H_t$ starting from $m = \delta_{(\bar{x},\bar{y})}$. The sides of the obtuse triangle are not perpendicular to each other; this and the fact that $w(s)$ and $z(s)$ are "typical" Brownian paths ensures that $K(s)$ is defined in a unique way for all $s$ a.s. and that we will never have $w(s) = z(s)$. The direction of $K$ either remains constant or approaches the direction of the side which is currently reflecting one of the paths [BB97].

For all $s$, the angle $\angle K(s)$ can never leave the interval $[\max(b,d) - \frac{\pi}{2}, \min(a,c) + \frac{\pi}{2}]$. We always have

$$w(s)^1 < z(s)^1, \text{ for all } s. \tag{9}$$

Let $\phi^1(x,y) = \phi(x)$ and $\phi^2(x,y) = \phi(y)$. Since every path in the support of $H$ starts at $(x,y)$ and $\phi$ satisfies the hypothesis of the theorem, by (9) we know that $\phi^1(w(s), z(s)) < \phi^2(w(s), z(s))$, for all $s$.

Using the relationship established earlier between super-reflected Brownian motions and solutions to (1), we may deduce

$$\begin{aligned}
u(t,x) &= -\log E(\exp(-\langle X_t, \phi \rangle)) \\
&= -\log E(\exp(-\int H_t^1(dy_1)\phi(y_1(t)))) \\
&= -\log E(\exp(-\int H_t(dy_1, dy_2)\phi^1(y_1(t), y_2(t)))) \\
&\leq -\log E(\exp(-\int H_t(dy_1, dy_2)\phi^2(y_1(t), y_2(t)))) \\
&= -\log E(\exp(-\int H_t^2(dy_2)\phi(y_2(t)))) \\
&= -\log E(\exp(-\langle Y_t, \phi \rangle)) = u(t,y).
\end{aligned} \tag{10}$$

Hence the solution $u(t,x)$ is monotonically increasing in $x^1$ for all $(x^1, x^2) \in K \cap D$. Since this is true for every line $K$ satisfying (8), the gradient of $u(t,x)$ must satisfy the condition stated in the theorem. $\square$

**Remarks:** The assumption that $D$ is a triangle plays no role in the arguments described in the proof of Theorem 1.1. The only property $D$ needs to satisfy is that, if two reflected Brownian motions in the domain are synchronous coupled, then the left particle will stay left of the other particle for all time $t$. We refer the reader to [BB97] or [Ath98] for examples of certain polygonal and non-convex domains $D$ in $\mathbb{R}^2$ with the above property.

We also wish to point out that the special form of $\Phi$, which enabled the particle representation of the partial differential equation (1), was crucial for the result to hold true.

The Feynman-Kacs representation of solutions to the partial differential equation do not yield monotonicity properties of the solution. Hence there does not seem to be an obvious probabilistic method to extend the above results to other partial differential equations.

**Acknowledgments** : I would like to thank my thesis advisor, Prof. K. Burdzy for suggesting this problem to me and for all the discussions that went into resolving the technical details.